\documentclass[10pt]{amsart}
\usepackage{latexsym}   
\usepackage{amssymb}    
\usepackage{amsmath}    
\usepackage{amsthm}     
\usepackage[all]{xy}
\usepackage{hyperref}
\newcommand{\pa}{\partial}

\newcommand{\del}{\delta}

\newcommand{\la}{\lambda}\newcommand{\La}{\Lambda}\newcommand{\om}{\omega}

\renewcommand{\thefootnote}

\newtheorem{theorem}{Theorem}[section]

\theoremstyle{definition}

\theoremstyle{remark}

\numberwithin{equation}{section}

\title[The B\"{a}cklund transforms of Peterson's isometric deformations of diagonal higher dimensional quadrics without center] {The B\"{a}cklund transforms of Peterson's isometric deformations of diagonal higher dimensional quadrics without center}

\author[  Ion I. Dinc\u{a}]{Ion I. Dinc\u{a}}
\address{Department of Applied Mathematics, Faculty of Applied Sciences,
National University of Science and Technology Politehnica Bucharest 313 Spl. Independentei 060042
Bucharest, Romania}
 \email{idinca@upb.ro}

\begin{document}


\begin{abstract}
We provide the B\"{a}cklund transforms of Peterson's isometric deformations of diagonal higher dimensional quadrics without center. These are found explicitly. can be iterated via the Bianchi Permutability Theorem and can be further iterated via the $3$-moving M\"{o}bius configuration, thus producing explicit solutions depending on arbitrarily many constants.
\end{abstract}

\maketitle

\tableofcontents \pagenumbering{arabic}

\section{Introduction}

The classical problem of finding the isometric deformations of surfaces (see Eisenhart \cite{E1})
was stated in 1859 by the French Academy of Sciences as

\begin{center}
{\it To find all surfaces applicable to a given one.}
\end{center}

Probably the most successful researcher of this problem is Bianchi, who in 1906 in \cite{B1}
solved the problem for quadrics by introducing the B\"{a}cklund transformation of surfaces
isometric to quadrics and the isometric correspondence provided by the Ivory affine transformation.
By 1909 Bianchi had a fairly complete treatment \cite{B2}; see also Eisenhart  \cite{E2} for transformation of surfaces.

In 1919-1920 Cartan has shown in \cite{C} (using mostly projective arguments and his exterior
differential systems in involution and exteriorly orthogonal forms tools) that space forms of
dimension $n$ admit rich families of isometric deformations (depending on $n(n-1)$ functions of one
variable) in surrounding space forms of dimension $2n-1$, that such isometric deformations admit
lines of curvature (given by a canonical form of exteriorly orthogonal forms; thus they have flat
normal bundle; since the lines of curvature on $n$-dimensional space forms (when they are considered
by definition as quadrics in surrounding $(n+1)$-dimensional space forms) are undetermined, the
lines of curvature on the isometric deformation and their corresponding curves on the quadric provide
the common conjugate system) and that the codimension $n-1$ cannot be lowered without obtaining
rigidity as the isometric deformation being the defining quadric.

In 1979, upon a suggestion from S. S. Chern and using Chebyshev coordinates on
$\mathbb{H}^n(\mathbb{R})$ (the Cartan-Moore Theorem; they are lines of curvature and thus in
bijective correspondence with isometric deformations of $\mathbb{H}^n(\mathbb{R})$ in
$\mathbb{R}^{2n-1}$) Tenenblat-Terng have developed in \cite{TT} the B\"{a}cklund transformation of
$\mathbb{H}^n(\mathbb{R})$ in $\mathbb{R}^{2n-1}$ (and Terng in \cite{T} has developed the {\it
Bianchi Permutability Theorem} for this B\"{a}cklund transformation).

In 1983 Berger, Bryant and Griffiths \cite{BBG} proved (including by use of tools from algebraic
geometry) in particular the Cartan's essentially projective arguments (including the exterior part
of his exteriorly orthogonal forms tool) can be used to generalize his results to $n$-dimensional
general quadrics with positive definite metric (thus they can appear as quadrics in
$\mathbb{R}^{n+1}$ or as space-like quadrics in $\mathbb{R}^n\times(i\mathbb{R})$) admitting rich
families of isometric deformations (depending on $n(n-1)$ functions of one variable) in surrounding
Euclidean space $\mathbb{R}^{2n-1}$, that the codimension $n-1$ cannot be lowered without obtaining
rigidity as the isometric deformation being the defining quadric and that quadrics are the only
Riemannian $n$-dimensional manifolds that admit a family of isometric deformations in
$\mathbb{R}^{2n-1}$ as rich as possible for which the exteriorly orthogonal forms tool (naturally
appearing from the Gau\ss\ equations) can be applied.

All computations are local and assumed to be valid on their open domain of validity without further
details; all functions have the assumed order of differentiability and are assumed to be invertible,
non-zero, etc when required (for all practical purposes we can assume all functions to be analytic).

Consider the complexified Euclidean space
$$(\mathbb{C}^p,<.,.>),\ <x,y>:=x^Ty,\ |x|^2:=x^Tx,\ x,y\in\mathbb{C}^p$$
with standard basis $\{e_j\}_{j=1,...,p},\ e_j^Te_k=\del_{jk}$.

Isotropic (null) vectors are those vectors $v$ of length $0$ ($|v|^2=0$); since most vectors are not
isotropic we shall call a vector simply vector and we shall only emphasize isotropic when the vector
is assumed to be isotropic. The same denomination will apply in other settings: for example we call
quadric a non-degenerate quadric (a quadric projectively equivalent to the complex unit sphere).

A quadric $x\subset\mathbb{C}^{n+1}$ is given by the quadratic equation
$$Q(x):=\begin{bmatrix}x\\1\end{bmatrix}^T\begin{bmatrix}A&B\\B^T&C\end{bmatrix}
\begin{bmatrix}x\\1\end{bmatrix}=x^T(Ax+2B)+C=0,$$
$$A=A^T\in\mathbf{M}_{n+1}(\mathbb{C}),\ B\in\mathbb{C}^{n+1},\ C\in\mathbb{C},\
\begin{vmatrix}A&B\\B^T&C\end{vmatrix}\neq 0.$$

The quadric $x\subset\mathbb{C}^{n+1}$ is degenerate if $\begin{vmatrix}A&B\\B^T&C\end{vmatrix}=0$.

There are many definitions of totally real (sub)spaces of $\mathbb{C}^{n+1}$, some even involving a
hermitian inner product, but all definitions coincide: an $(n+1)$-totally real subspace of
$\mathbb{C}^{n+1}$ is of the form $(R,t)(\mathbb{R}^k\times(i\mathbb{R})^{n+1-k}),\ k=0,...,n+1$,
where $(R,t)\in\mathbf{O}_{n+1}(\mathbb{C})\ltimes\mathbb{C}^{n+1}$. Now a totally real quadric is
simply an $n$-dimensional quadric in an $(n+1)$-totally real subspace of $\mathbb{C}^{n+1}$.

A metric classification of all (totally real) quadrics in $\mathbb{C}^{n+1}$ requires the notion of
symmetric Jordan canonical form of a symmetric complex matrix. The symmetric Jordan blocks
are:
$$J_1:=0=0_{1,1}\in\mathbf{M}_1{\mathbb{C}},\ J_2:=f_1f_1^T\in\mathbf{M}_2{\mathbb{C}},\
J_3:=f_1e_3^T+e_3f_1^T\in\mathbf{M}_3{\mathbb{C}},$$
$$J_4:=f_1\bar f_2^T+f_2f_2^T+\bar f_2f_1^T\in\mathbf{M}_4{\mathbb{C}},\
J_5:=f_1\bar f_2^T+f_2e_5^T+e_5f_2^T+\bar f_2f_1^T\in\mathbf{M}_5{\mathbb{C}},$$
$$J_6:=f_1\bar f_2^T+f_2\bar f_3^T+f_3f_3^T+\bar f_3f_2^T+\bar f_2f_1^T\in\mathbf{M}_6{\mathbb{C}},$$
etc,where $f_j:=\frac{e_{2j-1}-ie_{2j}}{\sqrt{2}}$ are the standard isotropic vectors (at least the
blocks $J_2,\ J_3$ were known to the classical geometers).

Any symmetric complex matrix can be brought via conjugation with a complex rotation to the symmetric
Jordan canonical form, that is a matrix block decomposition with blocks of the form $a_jI_p+J_p$;
totally real quadrics are obtained for eigenvalues $a_j$ of the quadratic part $A$ defining the
quadric being real or coming in complex conjugate pairs $a_j,\ \bar a_j$ with subjacent symmetric
Jordan blocks of same dimension $p$. Just as the usual Jordan block $\sum_{j=1}^pe_je_{j+1}^T$ is
nilpotent with $e_{p+1}$ cyclic vector of order $p$, $J_p$ is nilpotent with $\bar f_1$ cyclic
vector of order $p$, so we can take square roots of symmetric Jordan matrices without isotropic kernels
($\sqrt{aI_p+J_p}:=\sqrt{a}\sum_{j=0}^{p-1}(_j^{\frac{1}{2}})a^{-j}J_p^j,\ a\in\mathbb{C}^*,\
\sqrt{a}:=\sqrt{r}e^{i\theta}$ for $a=re^{2i\theta},\ 0<r,\ -\pi<2\theta\le\pi$), two matrices with
same symmetric Jordan decomposition type (that is $J_p$ is replaced with a polynomial in $J_p$) commute, etc.

The confocal family $\{x_z\}_{z\in\mathbb{C}}$ of a quadric $x_0\subset\mathbb{C}^{n+1}$ in
canonical form (depending on as few constants as possible) is given in the projective space
$\mathbb{CP}^{n+1}$ by the equation
$$Q_z(x_z):=\begin{bmatrix}x_z\\1\end{bmatrix}^T(\begin{bmatrix}A&B\\B^T&C\end{bmatrix}^{-1}-
z\begin{bmatrix}I_{n+1}&0\\0^T&0\end{bmatrix})^{-1}\begin{bmatrix}x_z\\1\end{bmatrix}=0,$$
where

$\bullet\ A=A^T\in\mathbf{GL}_{n+1}(\mathbb{C})\ \mathrm{symmetric\ Jordan},\ B=0\in\mathbb{C}^{n+1},\ C=-1$ for
quadrics with center.

$\bullet\ A=A^T\in\mathbf{M}_{n+1}(\mathbb{C})\ \mathrm{symmetric\ Jordan},\ \ker(A)=\mathbb{C}e_{n+1},\
B=-e_{n+1},\ C=0$ for quadrics without center and

$\bullet\ A=A^T\in\mathbf{M}_{n+1}(\mathbb{C})\ \mathrm{symmetric\ Jordan},\ \ker(A)=\mathbb{C}f_1,\ B=-\bar f_1,\
C=0$ for isotropic quadrics without center.

From the definition one can see that the family of quadrics confocal to $x_0$ is the adjugate of the
pencil generated by the adjugate of $x_0$ and Cayley's absolute $C(\infty)\subset\mathbb{CP}^n$ in
the hyperplane at infinity and it depends on the spectral parameter $z$; since Cayley's absolute encodes the Euclidean structure of
$\mathbb{C}^{n+1}$ (it is the set invariant under rigid motions and homotheties of
$\mathbb{C}^{n+1}:=\mathbb{CP}^{n+1}\setminus\mathbb{CP}^n$) the mixed {\it metric-projective}
character of the confocal family becomes clear.

For quadrics with center $\mathrm{spec}(A)$ is unambiguous (does not changes under rigid motions $(R,t)\in
\mathbf{O}_{n+1}(\mathbb{C})\ltimes\mathbb{C}^{n+1}$) but for (isotropic) quadrics without center it may change with
$(p+1)$-roots of unity for the block of ($f_1$ in $A$ being $J_p$) $e_{n+1}$ in $A$ being $J_1$ even
under rigid motions which preserve the canonical form, so it is unambiguous up to $(p+1)$-roots of
unity; for simplicity we make a choice and work with it.

We have the diagonal quadrics with(out) center respectively for $A=\sum_{j=1}^{n+1}a_j^{-1}e_je_j^T,\
A=\sum_{j=1}^na_j^{-1}e_je_j^T$; the diagonal isotropic quadrics without center come in different flavors, according to the
block of $f_1:\ A=J_p+\sum_{j=p+1}^{n+1}a_j^{-1}e_je_j^T$; in particular if $A=J_{n+1}$, then
$\mathrm{spec}(A)=\{0\}$ is unambiguous. General quadrics are those for which all eigenvalues have
geometric multiplicity $1$; equivalently each eigenvalue has an only corresponding symmetric Jordan block; in this
case the quadric also admits elliptic coordinates.

There are continuous groups of symmetries which preserve the symmetric Jordan canonical form for more than one symmetric Jordan
block corresponding to an eigenvalue, so from a metric point of view a metric classification
according to the elliptic coordinates and continuous symmetries may be a better one.

With $R_z:=I_{n+1}-zA,\ z\in\mathbb{C}\setminus\mathrm{spec}(A)^{-1}$ the family of quadrics
$\{x_z\}_z$ confocal to $x_0$ is given by $Q_z(x_z)=x_z^TAR_z^{-1}x_z+2(R_z^{-1}B)x_z+C+
zB^TR_z^{-1}B=0$. For $z\in\mathrm{spec}(A)^{-1}$ we obtain singular confocal quadrics; those with
$z^{-1}$ having geometric multiplicity $1$ admit a singular set which is an $(n-1)$-dimensional
quadric projectively equivalent to $C(\infty)$, so they will play an important r\^{o}le in the
discussion of homographies $H\in\mathbf{PGL}_{n+1}(\mathbb{C})$ taking a confocal family into
another one, since $H^{-1}(C(\infty)),\ C(\infty)$ respectively $C(\infty),\ H(C(\infty))$ will
suffice to determine each confocal family. Such homographies preserve all metric-projective
properties of confocal quadrics (including the good metric properties of the Ivory affine
transformation) and thus all integrable systems whose integrability depends only on the family of
confocal quadrics (the resulting involutory transformation between integrable systems is called
{\it Hazzidakis} (H) by Bianchi for the integrable system in discussion being the problem of
isometrically deforming quadrics).

The Ivory affine transformation is an affine correspondence between confocal quadrics and having
good metric properties: it is given by
$$x_z=\sqrt{R_z}x_0+C(z),\ C(z):=-(\frac{1}{2}\int_0^z(\sqrt{R_w})^{-1}dw)B.$$
Note that $C(z)=0$ for quadrics with center, $=\frac{z}{2}e_{n+1}$ for quadrics without center; for isotropic quadrics without center it is the Taylor series of
$\frac{1}{2}\int_0^z(\sqrt{1-w})^{-1}dw$ at $z=0$ with each monomial $z^{k+1}$ replaced by
$z^{k+1}J_p^k\bar f_1$, where $J_p$ is the block of $f_1$ in $A$ and thus a polynomial of degree $p$
in $z$.

Peterson's isometric deformations in $\mathbb{C}^{2n-1}$ of an $n$-dimensional quadric $x_0\subset\mathbb{C}^{n+1}$
is a maximal $(n-1)$-dimensional explicit family of non-trivial isometric deformations of $x_0$ with common conjugate system and non-degenerate joined second fundamental forms (see Dinc\u{a} \cite{D3}).

In Dinc\u{a} \cite{D2} we have found the B\"{a}cklund transforms of Peterson's isometric deformations of $2$-dimensional diagonal quadrics without center; the purpose of this note is to generalize this result to higher dimensions.

For the specific computations of isometric deformations of diagonal quadrics without center we shall use the convention
$\mathbb{C}^n\subset\mathbb{C}^{n+1}$ with $0$ on the $(n+1)^{\mathrm{th}}$ component; thus for
example we can multiply $(n+1,n+1)$-matrices with $n$-column vectors and similarly one can extend
$(n,n)$-matrices to $(n+1,n+1)$-matrices with zeroes on the last column and row. The converse is
also valid: an $(n+1,n+1)$-matrix with zeroes on the last column and row (or multiplied on the left
with an $n$-row vector and on the right with an $n$-column vector) will be considered as an
$(n,n)$-matrix.

\subsection{The main theorems from Dinc\u{a} \cite{D1}}\ \

With $V:=\sum_{k=1}^nv^ke_k=[v^1\ ...\ v^n]^T$ consider the complex paraboloid
$$Z=Z(v^1,...,v^n)=V+\frac{|V|^2}{2}e_{n+1}.$$

We have the diagonal quadrics without center $x_0:=LZ,\ L\in\mathbf{GL}_{n+1}(\mathbb{C}),$
$$L:=(\sqrt{A+e_{n+1}e_{n+1}^T})^{-1},\ A:=\mathrm{diag}[a_1^{-1}. . . a_n^{-1}\ 0],\ a_j\in\mathbb{C}\setminus\{0\},\ j=1,...,n.$$

Consider the non-trivial isometric deformation $x\subset\mathbb{C}^{2n-1}$ of the diagonal quadric without center $x_0\subset\mathbb{C}^{n+1}$ with common conjugate system $(u^1,...,u^n)$.

Because $(v^1,...,v^n)$ are isothermal-conjugate and $(u^1,...,u^n)$ are conjugate on $x_0$, the
Jacobian $\frac{\pa(v^1,...,v^n)}{\pa(u^1,...,u^n)}$ has orthogonal columns, so with
$$\la_j:=|\pa_{u^j}V|\neq 0,\ \La:=[\la_1\ ...\ \la_n]^T,\ \del:=\mathrm{diag}[du^1\ ...\ du^n]$$
we have $dV=R\del\La,\ R\subset\mathbf{O}_n(\mathbb{C})$.

We have the next theorem from Dinc\u{a} \cite{D1} on non-trivial isometric deformations of diagonal higher dimensional quadrics without center:

\begin{theorem}\label{th:th4}
Non-trivial isometric deformations $x\subset\mathbb{C}^{2n-1}$ of the diagonal quadric without center $x_0\subset\mathbb{C}^{n+1}$
with isothermal-conjugate system $(v^1,...,v^n)$ and common conjugate system $(u^1,...,u^n)$ are
in bijective correspondence with solutions of the completely integrable linear differential system
(\ref{eq:dVqwc}) with compatibility condition (\ref{eq:domqwc}).
\end{theorem}

\begin{eqnarray}\label{eq:dVqwc}
dV=R\del\La,\ d\La=\om\La-\del R^TAV,\ \La^T\La=-(V^TAV+1),\nonumber\\
\om:=\sum_{j=1}^n(e_je_j^TR^T\pa_{u^j}R\del+\del R^T\pa_{u^j}Re_je_j^T).
\end{eqnarray}

\begin{eqnarray}\label{eq:domqwc}
d\om-\om\wedge\om=-\del R^TAR\wedge\del,\ \om\wedge\del-\del\wedge R^TdR=0\Leftrightarrow\nonumber\\
e_j^T[\pa_{u^j}(R^T\pa_{u^j}R)-\pa_{u^k}(R^T\pa_{u^k}R)-\sum_lR^T\pa_{u^l}Re_le_l^TR^T\pa_{u^l}R
+R^TAR]e_k=0,\ j\neq k,\nonumber\\
e_j^TR^T\pa_{u^l}Re_k=0\ \mathrm{for}\ j,k,l\ \mathrm{distinct},\ R\subset\mathbf{O}_n(\mathbb{C}).
\end{eqnarray}

We have the next three theorems from Dinc\u{a} \cite{D1} on the isometric deformations of higher dimensional quadrics which are the
perfect analogous of Bianchi's main three theorems on the isometric deformations of $2$-dimensional quadrics:

\begin{theorem}\label{th:th1}
{\bf I (existence and inversion of the B\"{a}cklund transformation for quadrics and the isometric
correspondence provided by the Ivory affine transformation)}

Any non-trivial isometric deformation $x^0\subset\mathbb{C}^{2n-1}$ of an $n$-dimensional sub-manifold
$x_0^0\subset x_0$ ($x_0\subset\mathbb{C}^{n+1}\subset\mathbb{C}^{2n-1}$ being a quadric) appears as
a focal sub-manifold of an $(\frac{n(n-1)}{2}+1)$-dimensional family of Weingarten congruences, whose
other focal sub-manifolds $x^1=B_z(x^0)$ are isometric, via the Ivory affine transformation between
confocal quadrics, to sub-manifolds $x_0^1$ in the same quadric $x_0$. The determination of these
sub-manifolds requires the integration of a family of Ricatti equations depending on the parameter
$z$ (we ignore for simplicity the dependence on the initial value data in the notation $B_z$).
Moreover, if we compose the inverse of the rigid motion provided by the Ivory affine transformation
$(R_0^1,t_0^1)$ with the rolling of $x_0^0$ on $x^0$, then we obtain the rolling of $x_0^1$
on $x^1$ and $x^0$ reveals itself as a $B_z$ transformation of $x^1$.
\end{theorem}

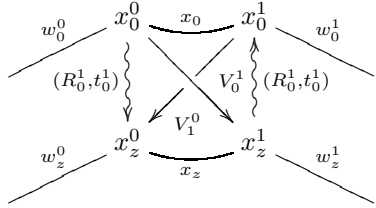
\begin{center}
$\xymatrix@!0{&&x_0^0\ar@{-}[drdr]\ar@/_/@{-}[rr]^{x_0}\ar@{~>}[dd]_{(R_0^1,t_0^1)}&&
x_0^1\ar@{<~}[dd]^{(R_0^1,t_0^1)}&&\\
\ar@{-}[urr]^{w_0^0}&&&\ar[dl]^>>>>{V_1^0}\ar[dr]^>>>>>{V_0^1}&&&\ar@{-}[ull]_{w_0^1}\\&&
x_z^0\ar@{-}'[ur][urur]\ar@/_/@{-}[rr]_{x_z}&&x_z^1&\\\ar@{-}[urr]^{w_z^0}&&&
&&&\ar@{-}[ull]_{w_z^1}}$
\end{center}

\begin{eqnarray}\label{eq:rjtj}
(R_j,t_j)(x_0^j,dx_0^j):=(R_jx_0^j+t_j,R_jdx_0^j)=(x^j,dx^j),\
(R_j,t_j)\subset\mathbf{O}_{2n-1}(\mathbb{C})\ltimes\mathbb{C}^{2n-1},\ j=0,1,\nonumber\\
(R_0^1,t_0^1)=(R_1,t_1)^{-1}(R_0,t_0).
\end{eqnarray}
\begin{theorem}\label{th:th2}
{\bf II (Bianchi Permutability Theorem)}

If $x^0\subset\mathbb{C}^{2n-1}$ is a non-trivial isometric deformation of an $n$-dimensional sub-manifold
$x_0^0\subset x_0$ ($x_0\subset\mathbb{C}^{n+1}\subset\mathbb{C}^{2n-1}$ being a(n isotropic) quadric without center) and $x^1=B_{z_1}(x^0),\ x^2=B_{z_2}(x^0)$, then one can find only by algebraic computations a
sub-manifold $B_{z_2}(x^1)=x^3=B_{z_1}(x^2)$; thus $B_{z_2}\circ B_{z_1}=B_{z_1}\circ B_{z_2}$ and
once all B transforms of the seed $x^0$ are found, the B transformation can be iterated using only
algebraic computations.
\end{theorem}

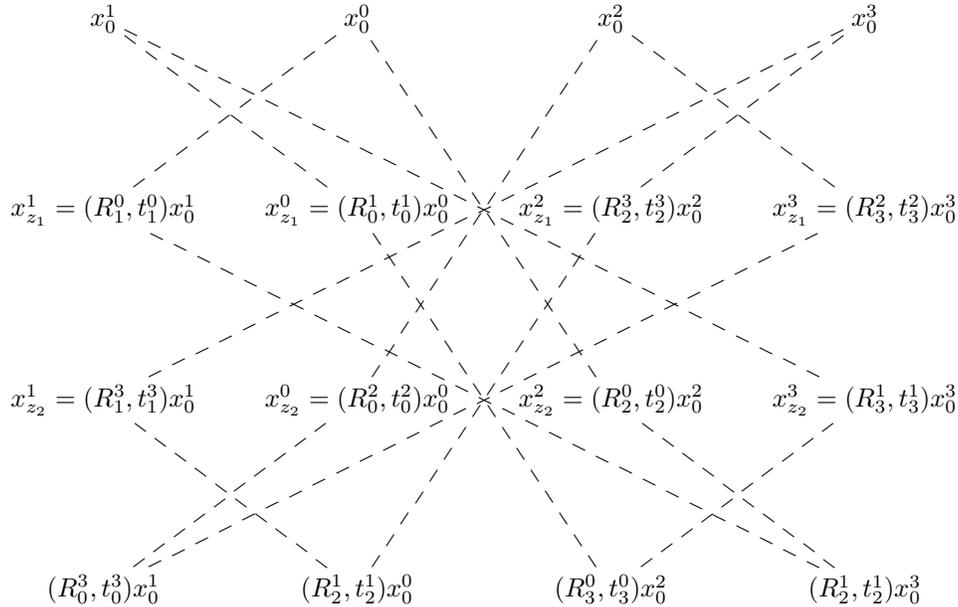
\begin{center}
$\xymatrix@!0{x_0^1\ar@{--}[dddrrrr]\ar@{--}[ddddddrrrrrrrrrrrr]&&&&
x_0^0\ar@{--}[dddllll]\ar@{--}[ddddddrrrr]&&&&
x_0^2\ar@{--}[ddddddllll]\ar@{--}[dddrrrr]&&&&
x_0^3\ar@{--}[ddddddllllllllllll]\ar@{--}[dddllll]\\
\\  \\
x_{z_1}^1=(R_1^0,t_1^0)x_0^1\ar@{--}[ddddddrrrrrrrrrrrr]&&&&
x_{z_1}^0=(R_0^1,t_0^1)x_0^0\ar@{--}[ddddddrrrr]&&&&
x_{z_1}^2=(R_2^3,t_2^3)x_0^2\ar@{--}[ddddddllll]&&&&
x_{z_1}^3=(R_3^2,t_3^2)x_0^3\ar@{--}[ddddddllllllllllll]\\
\\  \\
x_{z_2}^1=(R_1^3,t_1^3)x_0^1\ar@{--}[dddrrrr]&&&&
x_{z_2}^0=(R_0^2,t_0^2)x_0^0\ar@{--}[dddllll]&&&&
x_{z_2}^2=(R_2^0,t_2^0)x_0^2\ar@{--}[dddrrrr]&&&&
x_{z_2}^3=(R_3^1,t_3^1)x_0^3\ar@{--}[dddllll]\\
\\   \\
(R_0^3,t_0^3)x_0^1&&&&(R_2^1,t_2^1)x_0^0&&&&(R_3^0,t_3^0)x_0^2&&&&(R_2^1,t_2^1)x_0^3}$
\end{center}

\begin{eqnarray}\label{eq:rjtjbpt}
(R_j,t_j)(x_0^j,dx_0^j)=(x^j,dx^j),\ j=0,...,3,\nonumber\\
(R_j^k,t_j^k)=(R_k,t_k)^{-1}(R_j,t_j),\ (j,k)=(0,1),(0,2),(1,3),(2,3),\nonumber\\
(R_0^1,t_0^1)(R_2^0,t_2^0)=(R_3^0,t_3^0)=(R_3^1,t_3^1)(R_2^3,t_2^3),\
(R_1^0,t_1^0)(R_3^1,t_3^1)=(R_2^1,t_2^1)=(R_2^0,t_2^0)(R_3^2,t_3^2).\nonumber\\
\end{eqnarray}

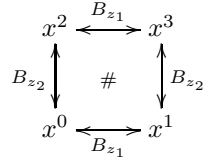
\begin{center}
$\xymatrix{\ar@{}[dr]|{\#}x^2\ar@{<->}[d]_{B_{z_2}}\ar@{<->}[r]^{B_{z_1}}&x^3\ar@{<->}[d]^{B_{z_2}}\\
x^0\ar@{<->}[r]_{B_{z_1}}&x^1}$
\end{center}

\begin{theorem}\label{th:th3}
{\bf III (existence of $3$-moving M\"{o}bius configurations $\mathcal{M}_3$)}

If $x^0\subset\mathbb{C}^{2n-1}$ is a non-trivial isometric deformation of an $n$-dimensional sub-manifold
$x_0^0\subset x_0$ ($x_0\subset\mathbb{C}^{n+1}\subset\mathbb{C}^{2n-1}$ being a(n isotropic) quadric without center) and $x^1=B_{z_1}(x^0),\ x^2=B_{z_2}(x^0),\ x^4=B_{z_3}(x^0)$ and by use of the Bianchi Permutability
Theorem one finds $B_{z_3}(x^2)=x^6=B_{z_2}(x^4),\ B_{z_1}(x^4)=x^5=B_{z_3}(x^1),\
B_{z_2}(x^1)=x^3=B_{z_1}(x^2),\ B_{z_3}(x^3)={x'}^7=B_{z_2}(x^5),\ B_{z_1}(x^6)={x''}^7=B_{z_3}(x^3),
\ B_{z_2}(x^5)={x'''}^7=B_{z_1}(x^6)$,then ${x'}^7={x''}^7={x'''}^7=:x^7$; thus once all $B$
transforms of the seed $x^0$ are found, the $B$ transformation can be further iterated using only
algebraic computations.
\end{theorem}

Existence of $3$-moving M\"{o}bius configuration imply existence of $p$-moving M\"{o}bius configurations, $p>3$.

\begin{center}
$\xymatrix@!0{&&&x^6\ar@{<->}[rrrr]^{B_{z_1}}&&&&x^7\\
&&&&&&&\\
x^4\ar@{<->}[uurrr]^{B_{z_2}} \ar@{<->}[rrrr]_{B_{z_1}}&&&&
x^5\ar@{<->}[uurrr]_>>>>>>>>>{B_{z_2}}&&&\\
&&&&&&&\\
&&&x^2\ar@{<->}'[r][rrrr]_{B_{z_1}}
\ar@{<->}'[uu][uuuu]_<<<<<{B_{z_3}}&&&&
x^3\ar@{<->}[uuuu]_{B_{z_3}}\\
&&&&&&&\\
x^0\ar@{<->}[rrrr]_{B_{z_1}}\ar@{<->}[uurrr]^{B_{z_2}}
\ar@{<->}[uuuu]^{B_{z_3}}&&&& x^1\ar@{<->}[uurrr]_{B_{z_2}}
\ar@{<->}[uuuu]_<<<<<<<<<<<<<<<<<<<<<<{B_{z_3}}&&&}$
\end{center}

Recall that the Ivory affine transformation between confocal diagonal quadrics without center is given by

$$x_z=\sqrt{R_z}x_0+\frac{z}{2}e_{n+1}=\sqrt{R_z}LZ+\frac{z}{2}e_{n+1}.$$

Consider two points $x_0^0,\ x_0^1\in x_0$ such that $x_0^0,\ x_z^1$ are in the symmetric tangency configuration

\begin{eqnarray}\label{eq:n00}
0=(N_0^0)^T(x_z^1-x_0^0)\Leftrightarrow|\sqrt{R_z}V_1-V_0|^2=-z(V_1^TAV_1+1)
\Leftrightarrow\nonumber\\
V_1^T\sqrt{R_z}V_0-\frac{|V_1|^2+|V_0|^2}{2}-\frac{z}{2}=0
\Leftrightarrow\nonumber\\
x_z^1=x_0^0+[\pa_{v_0^1}x_0^0\ ...\ \pa_{v_0^n}x_0^0](\sqrt{R_z}V_1-V_0).
\nonumber\\
\end{eqnarray}

Thus among the $2n$ functionally independent variables
$\{v_0^j,v_1^j\}_{j=1,...,n}$ a quadratic functional relation is established and only $2n-1$ among
them remain functionally independent:

\begin{eqnarray}\label{eq:dV0}
dV_0^T(\sqrt{R_z}V_1-V_0)=-dV_1^T(\sqrt{R_z}V_0-V_1).
\end{eqnarray}
Given a non-trivial isometric deformation $x^0\subset\mathbb{C}^{2n-1}$ of $x_0^0$ (that is $|dx^0|^2=
|dx_0^0|^2$) with orthonormal normal frame $N^0:=[N_{n+1}^0\ ...\ N_{2n-1}^0]$ consider the
$n$-dimensional sub-manifold

\begin{eqnarray}\label{x1x0}
x^1=x^0+[\pa_{v_0^1}x^0\ ...\ \pa_{v_0^n}x^0](\sqrt{R_z}V_1-V_0)\subset
\mathbb{C}^{2n-1}
\end{eqnarray}
(that is we restrict $\{v_1^j\}_{j=1,...,n}$ to depend only on the functionally independent
$\{v_0^j\}_{j=1,...,n}$ and constants in a manner that will subsequently become clear when we shall
impose the isometric correspondence provided by the Ivory affine transformation).

We take advantage now of the conjugate system$(u^1,...,u^n)$ common to $x_0^0,\ x^0$ and of the
non-degenerate joined second fundamental forms property; according to the principle of symmetry
$0\leftrightarrow 1$ we would like $(u^1,...,u^n)$ to be conjugate system to both $x^1$ and $x_0^1$
and also that the non-degenerate joined second fundamental forms property holds for $x_0^1,\ x^1$.
We obtain $dV_1=R_1\del\La_1$,

$$\La_1:=-\frac{1}{\sqrt{z}}R_0^T(\sqrt{R_z}V_1-V_0)$$

(note that the prime integral property $|\La_1|^2=-(V_1^TAV_1+1)$ is satisfied).

Thus

\begin{eqnarray}\label{eq:zr1}
\sqrt{z}R_1\La_0=\sqrt{R_z}V_0-V_1,\nonumber\\
-\sqrt{z}R_0\La_1=\sqrt{R_z}V_1-V_0
\end{eqnarray}
(the symmetry $(0,\sqrt{z})\leftrightarrow(1,-\sqrt{z})$ is required by (\ref{eq:dV0})). From the
first relation of (\ref{eq:zr1}) we get $V_1$ as an algebraic expression of $R_1,\ V_0,\ \La_0$;
replacing this into the second relation of (\ref{eq:zr1}) we get $\La_1$ as an algebraic expression
of $R_0,\ R_1,\ V_0,\ \La_0$:

\begin{eqnarray}\label{eq:V1La1}
V_1=\sqrt{R_z}V_0-\sqrt{z}R_1\La_0,\nonumber\\
\La_1=R_0^T(\sqrt{z}AV_0+\sqrt{R_z}R_1\La_0).
\end{eqnarray}
Differentiating the first relation of (\ref{eq:V1La1}) and using the second we get
$$R_1\del R_0^T(\sqrt{z}AV_0+\sqrt{R_z}R_1\La_0)=
\sqrt{R_z}R_0\del\La_0-\sqrt{z}[dR_1\La_0+R_1(\om_0\La_0-\del R_0^TAV_0)],$$
or the Ricatti equation

\begin{eqnarray}\label{eq:dR1}
dR_1=R_1\om_0+R_1\del R_0^TDR_1-DR_0\del,\ D:=\frac{\sqrt{R_z}}{\sqrt{z}}
\end{eqnarray}
in $R_1$.

We have now the theorem partially corresponding to {\bf Theorem} \ref{th:th1} for diagonal higher dimensional quadrics without center:

\begin{theorem}\label{th:th6}
Given $R_0\subset\mathbf{O}_n(\mathbb{C})$ solution of (\ref{eq:domqwc}) the Ricatti equation
(\ref{eq:dR1}) is completely integrable. If $R_1\subset\mathbf{M}_n(\mathbb{C})$ is a solution of
(\ref{eq:dR1}), then $R_1$ remains orthogonal if initially it was orthogonal and in this case it is
another solution of (\ref{eq:domqwc}) (thus producing an $(\frac{n(n-1)}{2}+1)$-dimensional family of
solutions). Moreover if $V_0,\ \La_0$ are solutions of (\ref{eq:dVqwc}) associated to $R_0$ (thus
producing a seed isometric deformation $x^0\subset\mathbb{C}^{2n-1}$ of $x_0^0$), then $V_1,\ \La_1$
given by

\begin{eqnarray}\label{eq:v1la1}
\begin{bmatrix}V_1\\\La_1\end{bmatrix}=\sqrt{z}\begin{bmatrix}I_n&0\\0&R_0^T\end{bmatrix}
\begin{bmatrix}D&-I_n\\A&D\end{bmatrix}\begin{bmatrix}I_n&0\\0&R_1\end{bmatrix}
\begin{bmatrix}V_0\\\La_0\end{bmatrix},\
D:=\frac{\sqrt{R_z}}{\sqrt{z}}
\end{eqnarray}
are solutions of (\ref{eq:dVqwc}) associated to $R_1$ (thus producing a leaf isometric deformation
$x^1\subset\mathbb{C}^{2n-1}$ of $x_0^1$) and we have the symmetry $(0,\sqrt{z})\leftrightarrow
(1,-\sqrt{z})$.
\end{theorem}

\subsection{The main theorem}\ \

Here we have the main Theorem of this paper:

\begin{theorem}\label{th:th7}
Peterson's isometric deformations of a diagonal higher dimensional  quadric without center can be found explicitly (modulo a quadratic matrix system). The B\"{a}cklund transforms of Peterson's isometric deformations of a diagonal higher dimensional  quadric without center can be found explicitly, can be iterated via the Bianchi Permutability Theorem and can be further iterated via the $3$-moving M\"{o}bius configuration, thus producing explicit solutions depending on arbitrarily many constants.
\end{theorem}

Note that by this main result one can also find the solitons of higher dimensional pseudo-spheres (see Tenenblat-Terng \cite{TT} and Terng \cite{T}).

In the next section we provide the proof of the main Theorem of this paper.

\section{Proof of the main  theorem}

\subsection{Some computations from \cite{D4}}\ \

First we reproduce some computations from \cite{D4}.

Note that once we know a solution
$R\subset\mathbf{O}_n(\mathbb{C})$ of (\ref{eq:domqwc}), a
solution $V,\ \La$ of (\ref{eq:dVqwc}) and certain linearly
independent solutions $V_j,\ \La_j,\ j=1,\cdots,2n-1$ of the
(homogeneous) differential part of (\ref{eq:dVqwc}), then we
can find up to rigid motions the space realization of the
isometric deformation $x\subset\mathbb{C}^{2n-1}$ of
$x_0$ up to quadratures and without the use of the
Gau\ss-Bonnet(-Peterson) Theorem (this is again due to Bianchi for
$n=2$).

Let $V_j,\ \La_j,\ j=1,\cdots,2n-1$ be certain linearly
independent solutions of the (homogeneous) differential part of
(\ref{eq:dVqwc}) and $V,\ \La$ solution of
(\ref{eq:dVqwc}) such that
$x\subset\mathbb{C}^{2n-1}$ given up to quadratures by
\begin{eqnarray}\label{eq:d'x}
dx:=[V_1\ \ \cdots\ \ V_{2n-1}]^TdV
\end{eqnarray}
is an isometric deformation of $x_0$ with non-degenerate joined
second fundamental forms. Note $d(dx)=0$. Since
$|dx|^2=|dx_0|^2,\ dx_0=(L+e_{n+1}V^T)dV$ we get
$$\sum_{j=1}^{2n-1}V_jV_j^T-VV^T=I_{1,n}L^TL.$$
Applying $d$ we get
$$\sum_{j=1}^{2n-1}V_j\La_j^T-V\La^T=0$$ (use $\del M=N\del,\
M,N\subset\mathbf{M}_n(\mathbb{C})\Leftrightarrow
M=N=\mathrm{diag}$). Applying $d$ again  we get
$$\sum_{j=1}^{2n-1}\La_j\La_j^T-\La\La^T=I_n.$$

With $\mathcal{V}:=[V_1\ \ \cdots\ \ V_{2n-1}\ \ iV],\
\mathcal{L}:=[\La_1\ \ \cdots\ \ \La_{2n-1}\ \ i\La]$ these can be
written as
\begin{eqnarray}\label{eq:vla0}
\begin{bmatrix}\mathcal{V}\\\mathcal{L}\end{bmatrix}[\mathcal{V}^T\
\ \mathcal{L}^T]=\begin{bmatrix}I_{1,n}L^TL& 0\\0&
I_n\end{bmatrix}.
\end{eqnarray}
We also have the prime integral property $[\mathcal{V}^T\ \
\mathcal{L}^T]\begin{bmatrix}A&0\\0&I_n\end{bmatrix}
\begin{bmatrix}\mathcal{V}\\\mathcal{L}\end{bmatrix}=
\begin{bmatrix}\mathcal{C}&ic\\ic^T&1\end{bmatrix},\
\mathcal{C}=\mathcal{C}^T\in\mathbf{M}_{2n-1}(\mathbb{C}),\\
c\in\mathbb{C}^{2n-1}$. Multiplying it on the left with
$\begin{bmatrix}\mathcal{V}\\\mathcal{L}\end{bmatrix}$ and using
(\ref{eq:vla0}) we get $\mathcal{C}=I_{2n-1},\ c=0$. Thus
$[\mathcal{V}^T\ \ \mathcal{L}^T]$ is determined modulo a
multiplication on the left with a rotation
$\begin{bmatrix}\mathcal{R}&0\\0&1\end{bmatrix}\in\mathbf{O}_{2n}(\mathbb{C})$
by
\begin{eqnarray}\label{eq:vla}
[\mathcal{V}^T\ \
\mathcal{L}^T]\begin{bmatrix}A&0\\0&I_n\end{bmatrix}
\begin{bmatrix}\mathcal{V}\\\mathcal{L}\end{bmatrix}=I_{2n}.
\end{eqnarray}
Thus  (\ref{eq:vla}) is equivalent to (\ref{eq:vla0}) and to
$[\mathcal{V}^TL^{-1}\ \
\mathcal{L}^T]\subset\mathbf{O}_{2n}(\mathbb{C})$.

For the non-degenerate joined second fundamental forms property we
need to prove that there is no vector field
$N\subset\mathbb{C}^{2n-1}$ along $x$ such that
$N^Tdx=0,\ -dN^Tdx=|dV|^2=\La^T\del^2\La$, that is the linear
system
$\begin{bmatrix}\mathcal{V}\\\mathcal{L}\end{bmatrix}\begin{bmatrix}N\\0\end{bmatrix}
=\begin{bmatrix}0\\\La\end{bmatrix}$ is inconsistent. Using
(\ref{eq:vla}) this becomes $\begin{bmatrix}N\\0\end{bmatrix}
=\mathcal{L}^T\La$, which is indeed inconsistent because
$|\La|^2\neq 0$.

\subsection{The main result}\ \

For Peterson's isometric deformation of diagonal higher dimensional quadrics without center we have $R=I_n,\ \om=0$ and (\ref{eq:domqwc}) is satisfied.

The relation (\ref{eq:dVqwc}) becomes

$$dV=\del\La,\ d\La=-\del AV,\ \La^T\La=-(V^TAV+1).$$

The solution of this differential system depends on $n-1$ constants (the prime integral removes a constant from the space of solutions and translations in $u^j,\ j=1,..,n$ another $n$ constants), so from the $2n$-dimensional space of the homogeneous part of the differential system we get the $(n-1)$- dimensionality of Peterson's isometric deformations of diagonal higher dimensional quadrics without center.

We have $\pa_{u^j}V=\la_je_j,\ \pa_{u^j}\La=-a_j^{-1}v^je_j$, so $v^j=v^j(u^j),\ \la_j=\la_j(u^j)$
and $\pa_{u^j}^2v^j+a_j^{-1}v^j=0$, so
$$v^j=c_1^je^{i\sqrt{a_j^{-1}}u^j}+c_2^je^{-i\sqrt{a_j^{-1}}u^j},$$
$$\la_j=i\sqrt{a_j^{-1}}(c_1^je^{i\sqrt{a_j^{-1}}u^j}-c_2^je^{-i\sqrt{a_j^{-1}}u^j})$$
and from the prime integral we get
$$\sum_{j=1}^n2a_j^{-1}(c_1^jc_2^j+c_2^jc_1^j)=-1.$$
Similarly with $V_l=[v_l^1\ ...\ v_l^n]^T,\ \La_l=[\la_{1l}\ ...\ \la_{nl}]^T,\ l=1,...,2n-1$ we have
$$v_l^j=c_{l1}^je^{i\sqrt{a_j^{-1}}u^j}+c_{l2}^je^{-i\sqrt{a_j^{-1}}u^j},$$
$$\la_{jl}=i\sqrt{a_j^{-1}}(c_{l1}^je^{i\sqrt{a_j^{-1}}u^j}-c_{l2}^je^{-i\sqrt{a_j^{-1}}u^j})$$ and the prime integral property (\ref{eq:vla}) becomes:

\begin{eqnarray}\label{eq:primi}
\sum_{j=1}^n2a_j^{-1}(c_{l1}^jc_{m2}^j+c_{l2}^jc_{m1}^j)=\del_{lm},\ l,m=1,...2n-1,\nonumber\\
\sum_{j=1}^n2a_j^{-1}(c_{l1}^jc_2^j+c_{l2}^jc_1^j)=0,\ l=1,...,2n-1,\nonumber\\
\sum_{j=1}^n2a_j^{-1}(c_1^jc_2^j+c_2^jc_1^j)=-1.
\end{eqnarray}

Thus Peterson's isometric deformations of a diagonal higher dimensional  quadric without center can be found explicitly modulo the quadratic matrix system (\ref{eq:primi}).

For $x=[x^1\ ...\ x^{2n-1}]$ Peterson's isometric deformations of a diagonal higher dimensional  quadric without center we have $$dx^l=V_l^TdV=\sum_{j=1}^nv_l^jdv^j=$$
$$=\sum_{j=1}^ni\sqrt{a_j^{-1}}(c_{l1}^jc_1^je^{2i\sqrt{a_j^{-1}}u^j}
-c_{l2}^jc_2^je^{-2i\sqrt{a_j^{-1}}u^j}+c_{l2}^jc_1^j-c_{l1}^jc_2^j)du^j,\ l=1,...,2n-1,$$  so
$$x^l=\sum_{j=1}^n[\frac{1}{2}(c_{l1}^jc_1^je^{2i\sqrt{a_j^{-1}}u^j}
+c_{l2}^jc_2^je^{-2i\sqrt{a_j^{-1}}u^j})+i\sqrt{a_j^{-1}}(c_{l2}^jc_1^j-c_{l1}^jc_2^j)u^j]+c^l,\ l=1,..., 2n-1.$$

From (\ref{eq:dR1}) with $R_0=I_n$ for Peterson's isometric deformations of a diagonal higher dimensional  quadric without center we get the Riccati equation in $R_1$ with constant coefficients

\begin{eqnarray}\label{eq:dR1p}
dR_1=R_1\del DR_1-D\del,\ D:=\frac{\sqrt{R_z}}{\sqrt{z}}
\end{eqnarray}
and $R_1$ remains orthogonal if initially it was orthogonal.

Note that (\ref{eq:dR1p}) admits the particular solution $R_1=I_n$, so for the solution $R_1:=I_n+U$ of (\ref{eq:dR1p}) $U$ will satisfy the matrix Bernoulli equation
$$dU=U\del D+\del DU+U\del DU.$$

With the substitution $Z:=U^{-1}$ we have the linear differential equation in $Z$:

$$-dZ=\del DZ+Z\del D+\del D.$$

With $Z=\{z_{jk}\}_{j,k=1,...,n}$ and $d_j:=\frac{\sqrt{1-za_j^{-1}}}{\sqrt{z}},\ j=1,..,n$
we have

$$-\pa_{u^j}z_{kl}=\del_{jk}d_jz_{jl}+z_{kj}d_j\del_{jl}+d_j\del_{jk}\del_{jl},\ j,k,l=1,...,n,$$
so
\begin{eqnarray}\label{eq:dZ}
\pa_{u^j}z_{jj}=-d_j(2z_{jj}+1),\ j=1,...,n,\nonumber\\
\pa_{u^j}z_{ll}=0,\ j,l=1,...,n\ j\neq l,\nonumber\\
\pa_{u^j}z_{jl}=-d_jz_{jl},\ j,l=1,...,n,\ j\neq l,\nonumber\\
\pa_{u^j}z_{lj}=-d_jz_{lj},\ j,l=1,...,n,\ j\neq l,\nonumber\\
\pa_{u^j}z_{kl}=0,\ j,k,l=1,...,n\ \mathrm{distinct}.
\end{eqnarray}

Once $Z$ can be found by quadratures from (\ref{eq:dZ}) the general solution of (\ref{eq:dR1p}) is given by
$R_1=I_n+Z^{-1}$; we get the B\"{a}cklund transformation for $R_1(0)\in\mathbf{O}_n(\mathbb{C})$.

From (\ref{eq:dZ}) we get $z_{jj}=-\frac{1}{2}+(z_{jj}(0)+\frac{1}{2})e^{-2d_ju^j},\ j=1,...,n,\
z_{jl}=z_{jl}(0)e^{-d_ju^j-d_lu^l},\ j,l=1,...,n,\ j\neq l$.

We have $R_1(0)\in\mathbf{O}_n(\mathbb{C})\Leftrightarrow\ Z(0)^T+Z(0)+I_n=0$, so $z_{jj}(0)=-\frac{1}{2},\ j=1,...,n,\ z_{jl}(0)=-z_{lj}(0),\ j,l=1,...,n,\ j\neq l.$

Note that the same argument for $D=\mathrm{diag}[\csc\sigma\ \ \cot\sigma\
... \ \cot\sigma]$ gives the solitons of higher dimensional pseudo-spheres, the vacuum soliton corresponding to $R_0=I_n,\ \om_0=0$ (see Tenenblat-Terng \cite{TT} and Terng \cite{T}; in their version $\om_0$ is replaced with $\om'_0$ and all computations are real). In Tenenblat-Terng's \cite{TT} version of (\ref{eq:dR1}) we need $R_0$ to have non-zero entries in the first row in order to get a good definition of $\om_0$, so a-priori (\ref{eq:dR1}) with $R_0=I_n,\ \om_0=0$ does not make much sense, but $R_1$ found from the integration of (\ref{eq:dR1p}) will mostly have non-zero entries on the first row, so the $1$-soliton $R_1$ is well defined. To see this by continuity we need $R_1(0)$ to mostly have non-zero entries in the first row; this condition is realized on the complement of a union of $n$ hyper-surfaces in $\mathbb{R}^{\frac{n(n-1)}{2}}$, so on an open dense set. To find the space realization of the solitons of higher dimensional pseudo-spheres we need to find the space realization of $0$ and $1$-solitons.

If $(R_1,z_1),(R_2,z_2)$ are solutions  of (\ref{eq:dR1p}), then by the Bianchi Permutability Theorem
we get (see Dinc\u{a} \cite{D1})

\begin{eqnarray}\label{eq:R3}
R_3=(D_2-D_1R_2R_1^T)(D_2R_2R_1^T-D_1)^{-1}
\end{eqnarray}
and $R_3$ is orthogonal.

If $(R_1,z_1),(R_2,z_2),(R_4,z_3)$ are solutions  of (\ref{eq:dR1p}), then by the existence of $3$-moving M\"{o}bius configurations with

$$\Box:=(D_2^2-D_3^2)D_1R_1+(D_3^2-D_1^2)D_2R_2+(D_1^2-D_2^2)D_3R_4$$

we get (see Dinc\u{a} \cite{D1})

$$D_1[(D_2^2-D_3^2)D_2R_3(D_2R_3-D_3R_5)^{-1}R_1-D_2^2R_1]=$$
$$(D_1^2D_2R_2-D_2^2D_1R_1)\Box^{-1}(D_2^2-D_3^2)(D_3R_4-D_1R_1)-D_2^2D_1R_1=$$
$$(D_1^2D_2R_2-D_2^2D_1R_1)\Box^{-1}(D_3^2-D_1^2)(D_2R_2-D_3R_3)-D_1^2D_2R_2=$$
$$D_2[(D_3^2-D_1^2)D_1R_3(D_3R_6-D_1R_3)^{-1}R_2-D_1^2R_2],$$
so the very left hand side and right hand side provide the good definition of and afford themselves the name
$D_1D_2D_3R_7$.

\end{document}